\documentclass [11pt,a4paper]{article}
\usepackage[utf8]{inputenc}

\usepackage[pdftex]{graphicx}
\usepackage{amsmath}
\usepackage{amsfonts}
\usepackage{amssymb}
\usepackage{graphics}
\usepackage{rotating}
\usepackage{pifont}
\usepackage{wasysym}
\usepackage{theorem}
\usepackage{url}
\usepackage{array}
\usepackage[arrow, curve, matrix]{xy}
\usepackage[top= 35mm, bottom=30mm, left=35mm, right=35mm]{geometry}
\usepackage{mathtools}
\usepackage{float}

\newtheorem{Satz}{Satz}

\newtheorem{Lemma}[Satz]{Lemma}
\newtheorem{Prop}[Satz]{Proposition}

\newtheorem{Sum}[Satz]{Summary}

\newtheorem{Char}[Satz]{Characterisation}

{\theorembodyfont{\rmfamily} \newtheorem{Def}[Satz]{Definition}
{\theorembodyfont{\rmfamily} \newtheorem{Rem}[Satz]{Remark}
{\theorembodyfont{\rmfamily} \newtheorem{Nota}[Satz]{Notation}

}

\DeclareMathOperator{\td}{trdeg}

\DeclareMathOperator{\llog}{llog}
\DeclareMathOperator{\lllog}{lllog}
\DeclareMathOperator{\Cond}{Co}
\DeclareMathOperator{\inter}{int}
\DeclareMathOperator{\Pole}{Pole}

\parindent 0pt
\parskip 3pt

\begin{document}
\renewcommand{\thefootnote}{\fnsymbol{footnote}}
\title{Geometric interpretations of a counterexample to Hilbert's 14th problem and rings of bounded polynomials on semialgebraic sets.}
\author{Sebastian Krug\footnotemark[1]} 
\date{}

\maketitle
\pagestyle{headings}

\newcommand{\mc}{\mathcal}
\newcommand{\GO}{\mathcal{O}}
\newcommand{\wh}{\widehat}
\newcommand{\m}{\mbox}
\newcommand{\mbb}{\mathbb}
\newcommand{\ul}{\underline}

\sloppy

\begin{abstract} 
We construct open semialgebraic subsets $S$ of $\mbb{R}^3$, such that the ring of bounded polynomials on $S$, written $B_{\mbb{R}^3}(S)$, is not finitely generated as an $\mbb{R}$-algebra. For this we use a family of counterexamples to Hilbert's 14th problem constructed by S. Kuroda: Let $R \subset K[x_1,x_2,x_3,x_4]$ be any algebra belonging to this family ($K$ any field of characteristic $0$). If $K= \mbb{R}$, we construct an open semialgebraic $S \subset \mbb{R}^3$ such that $B_{\mbb{R}^3}(S) \cong R$. Furthermore, for $K$ arbitrary, we construct an explicit smooth quasiprojective $K$-variety $U$, such that the ring of regular functions on $U$ is isomorphic to $R$.      
\end{abstract}

\footnotetext[1]{Universit\"at Hamburg, \texttt{email: sebastian.krug@uni-hamburg.de}}
\renewcommand{\thefootnote}{\arabic{footnote}}

Large parts of this article can be seen as a supplement to the article \cite{PS1} by D. Plaumann and C. Scheiderer, in which the authors
investigate the ring of bounded polynomials on a semialgebraic subset $S$ of an affine algebraic variety $V$ over the real numbers,   
$B_V(S):= \{ f \in \mbb{R}[V]\, | \, \text{$f_{|S}$ is bounded}\}$. Among other things, they introduce criteria for finite generatedness of $B_V(S)$ as an $\mbb{R}$-algebra. 

A semialgebraic  $S \subset V$ is by definition (see below) contained in the set of real points $V(\mbb{R})$ of $V$. We will say that a $K$-algebra is f.g. if it is finitely generated over the field $K$ (often $K= \mbb{R}$ understood). As shown in \cite{PS1}, for a semialgebraic $S$ which is neither Zariski-dense nor bounded, $B_V(S)$ is never f.g. 
On the other hand, $B_V(S)$ is always f.g., if $S$ is regular (e.g. open) in the euclidean topology of $V(\mbb{R})$ and $V$ nonsingular with $\dim \, V \le 2$.  
Here \emph{regular} means that $S$ is contained in the closure of its open interior, i.e. $S \subseteq \overline{\inter (S)}$. \footnote{$B_V(S)$ does not change when $S$ is altered inside a compact subset of the real points $V(\mbb{R})$. Hence, if one is interested in the dependence of $B_V(S)$ on geometric properties 
of $S$, one should define these properties in such a way, that they are independent of what happens inside such a compact subset. So the results just mentioned, would be better formulated using the adequate properties ``Zariski-dense \emph{at infinity}'' and ``regular \emph{at infinity}'', as done in \cite{PS1}. See Def. \ref{d1} (viii) below.}    
Even in dimension $>2$, regular semialgebraic sets $S$ such that $B_V(S)$ is not f.g. are not trivial to find,
but \cite{PS1} provides such examples, for instance regular semialgebraic subsets of certain smooth affine varieties $V$ of dimension $3$. 
However in these examples $V$ is never isomorphic to some $\mbb{A}_{\mbb{R}}^n$. 
So the question whether there is any regular semialgebraic $S \subset \mbb{R}^n$, such that $B_{\mbb{R}^n}(S)$ is not f.g., is not answered in \cite{PS1}.

The 14th problem of Hilbert is the question, whether for a field $K$, a polynomial ring $K[X_1,...,X_n]$ and an intermediate field $K|L|K(X_1,...,X_n)$, 
the intersection $R:= L \cap K[X_1,...,X_n]$ is always f.g. as a $K$-Algebra. The problem was inspired by questions about the finite generatedness of invariant subrings of 
$K[X_1,...,X_n]$ that Hilbert considered erroneously as already answered when he posed his 23 problems. The answer is ''no'' in general as shown by M. Nagata in \cite{MR0105409}. 
O. Zariski had shown before (\cite{MR0065217}), that for $\td_K L \le 2$, $R$ is always f.g., by interpreting the intersection $R$ as the ring of regular functions on a normal quasiprojective algebraic variety and showing that such rings are f.g. for varieties of dimension $\le 2$. 
(This theorem is applied in \cite{PS1}, to obtain the result for $S$ regular and $\dim V \le 2$ mentioned above.) 

For small $n >2$, and for $\td_K L = 3$, Hilbert's question remained open for a long time, until in 2004 S. Kuroda constructed a not f.g. example with $n=4$ 
and $\td_K L= 3$ (\cite{MR2078390}). 
Later he also produced an example with $n=3$ (\cite{MR2125538}). 

In the special case $K=\mbb{R}$, we interpret Kuroda's family of counterexamples from \cite{MR2078390} geometrically by constructing for each member $R$ of the family an  open semialgebraic $S \subset \mbb{R}^3$, such that $B_{\mbb{R}^3}(S)$ is isomorphic to the not f.g. $K$-algebra $R$. This answers the question left open in \cite{PS1}, mentioned above. We can choose $S$ to be basic open and defined by explicit inequalities. 
\footnote{This basic open semialgebraic set will be called $\widetilde{S}$ when we define it later in the article.}    
Furthermore, over any field $K$ of characteristic $0$ we explicitly construct quasiprojective varieties having Kuroda's not f.g. $K$-algebras as their rings of regular functions. They are obtained as open subvarieties of certain blowups of $\mbb{P}^3_K$. In the case $K=\mbb{R}$ the same construction provides completions of $\mbb{A}_{\mbb{R}}^3$ compatible with the semialgebraic subset $S$, in the sense of \cite{PS1} (see Def. \ref{d1} (ix) below).

In general there is an equivalence (see Appendix for details) between 
\begin{enumerate}
\item $K$-algebras $R$ of the form $R= A \cap L$ where $A$ is a normal f.g. $K$-algebra, $L|K$ a field extension, and 
\item $K$-algebras $R$ that appear as rings of global regular functions of irreducible normal (quasiprojective/quasiaffine) varieties. 
\end{enumerate}
This equivalence was discovered by O. Zariski and M. Nagata in their work on Hilbert's 14th problem. 
In the special case $K= \mbb{R}$, the class of these two equivalent types of $\mbb{R}$-algebras contains all $\mbb{R}$-algebras of a third type, namely rings of bounded polynomials of regular semialgebraic subsets of normal irreducible $\mbb{R}$-varieties. \footnote{Conversely, if some $R$ fulfills the two equivalent conditions, and has a (formally) real field of quotients then it is isomorphic to some such ring of bounded polynomials.} This is established (in dimension $\le 2$) by means of compatible completions in \cite{PS1}. In an Appendix we observe that it holds in arbitrary dimension. The main part of this article can be seen as an illustration of the connections between these three types of algebras at a concrete example.   

\textbf{Acknowledgments:} I would like to thank Claus Scheiderer, Daniel Plaumann and Emilie Dufresne for helpful comments and discussion. This article was written while I was employed at the
research training group ``analysis, geometry and string theory'' of the Leibniz Universit\"at Hannover.     

\section{Preliminaries}

We compile some definitions and results (mostly taken form \cite{PS1}) concerning $\mbb{R}$-varieties, semialgebraic sets and the ring of bounded polynomials on them. 

\begin{Def} \label{d1}
 
(i) For a field $K$, a \emph{$K$-variety} $V$ is a separated, reduced $K$-scheme of finite type. We write $\GO_V$ for the sheaf of regular $K$-valued functions on 
$V$, $\GO (V) := \GO_V(V)$  for the ring of global regular functions on $V$, and $K(V)$ for the field of rational functions on $V$. 
If $V$ is affine, then $K[V] = \GO (V)$ is the coordinate ring of $V$ over $K$. The set of non-singular points of $V$ is written $V_{reg}$. If $X$ is an affine (or projective) variety, $f_1,...,f_n$ elements of its (homogeneous) coordinate ring, then $\mc{V}_X(f_1,...,f_n)$ denotes the closed subvariety of $X$ defined by $f_1,...,f_n$. We sometimes call a subvariety $Y \subset X$ a \emph{divisor} if it is of pure codimension 1.  
 
(ii) For a $\mbb{R}$-variety $V$, $V(\mbb{R})$ is the \emph{set of $\mbb{R}$-rational points} of $V$. We always consider it endowed with the euclidean topology. 
 
(iii) An irreducible $\mbb{R}$-variety is said to be \emph{real}, if it has a nonsingular $\mbb{R}$-rational point, or equivalently, if $V(\mbb{R})$ 
is Zariski-dense in $V$, or equivalently if the field $\mbb{R}(X)$ is (formally) real.   

(iv) A \emph{semialgebraic subset} $S$ of an affine $\mbb{R}$-variety $V$ is a subset of $V(\mbb{R})$ of the form
\[S= \bigcup_{k=1}^n \{f_{k,1} >0,f_{k,2}>0, ..., f_{k,r_i}>0 , g_k = 0 \} \] 
where $n, r_i \in \mbb{N}$ and all $f_{k,j}, g_k \in \mbb{R}[V]$.  A semialgebraic subset of a general $\mbb{R}$-variety $V$ is an $S \subseteq V(\mbb{R})$ such that $S \cap W$ is semialgebraic in $W$ for each open affine $W \subseteq V$. 
 
(v) If $V$ is affine and $S \subset V$ is semialgebraic and open in $V(\mbb{R})$, one can omit the $g_k=0$, and if furthermore $n=1$ is possible, one calls the set $S$ \emph{basic open}. 
If we require $n=1$ and replace all the $>$ by $\ge$ we get what is called a \emph{basic closed} semialgebraic set. 
 
(vi) For $V$ an $\mbb{R}$-variety and $S \subset V$ semialgebraic, the \emph{ring of 
bounded polynomials} on $S$ (as a subset of $V$) is  $B_V(S):= \{ f \in \GO (V)\, | \, \text{$f_{|S}$ is bounded}\}$. \footnote{Maybe, if $V$ is not affine, one should call $B_V(S)$ ring of bounded regular functions instead.}      

(vii) A semialgebraic $S \subset V$ is called \emph{Zariski-dense at infinity}, if for every compact subset $C \subseteq V(\mbb{R})$, and for every closed algebraic subset $X$ of $V$:
\[S \subseteq C \cup X \; \Rightarrow \; X = V\]    

(viii) A semialgebraic $S \subset V$ is called \emph{regular at infinity}, if for some compact subset $C \subseteq V(\mbb{R})$:
\[ S \subseteq C \cup \overline{\inter(S \cap V_{reg}(\mbb{R}))} \quad \footnote{Where the overline denotes the closure in $V(\mbb{R})$, and $\inter$ the open interior. (Euclidian topology)}\] 
(Every non-bounded set that is regular at infinity is also Zariski-dense at infinity.)

(ix) For $V$ an irreducible and $X$ a complete $\mbb{R}$-variety, an open dense embedding  $V \hookrightarrow X$, is called a
\emph{completion} of $V$. For a semialgebraic  $S \subset V$, the completion is called \emph{compatible with $S$}, if 
for every irreducible component $Z$ of $X \smallsetminus V$:

\begin{enumerate}
\item $\GO_{X,Z}$ is a discrete valuation ring.  
\item If $\overline{S}$ is the closure of $S$ in $X(\mbb{R})$, then either $Z \cap \overline{S} = \emptyset$ or $Z \cap \overline{S}$ is Zariski-dense in $Z$.
\end{enumerate}
 
\end{Def}

\begin{Sum} \label{su} (Plaumann, Scheiderer)

For an irreducible $\mbb{R}$-variety $V$ and $S \subset V$ semialgebraic: 

(i) (Cor. 5.8. \cite{PS1}) If $V$ is affine, and $S$ is neither Zariski-dense at infinity nor bounded, then $B_V(S)$ is not f.g. over $\mbb{R}$ (and even not noetherian). 
 
(ii)(Thm. 3.8. \cite{PS1}) Let $V$ be normal. If $V \hookrightarrow X$ is a completion of $V$, compatible with $S$,
let $Z$ be the union of all irreducible components of $X \smallsetminus V$ whose intersection with $\overline{S}$ 
(cf. Definition \ref{d1} (ix)) is empty, and set $U:= X \smallsetminus Z$. Then, considering $\GO(U)$ in the natural way as a subalgebra of $\GO(V)$, one has
\[\GO(U) = B_V(S). \quad \footnote{The result (ii) is formulated in $\cite{PS1}$ only for affine $V$. But the proof given there works for general $V$ as well. (We will use the general case in the Appendix.)} \]

(iii)(Thm. 4.5. and 5.12. \cite{PS1}) If $V$ is normal and affine, $\dim V \le 2$, and $S$ is regular at infinity, then $V$ has a completion $V \hookrightarrow X$ compatible with $S$, and one can use this to show 
that $B_V(S)$ is f.g. over $\mbb{R}$.
\end{Sum}

\begin{Nota} (i) If $K[x_1,...,x_n]$ is some $K$-algebra, we often use the shorthand $K[\ul{x}]:= K[x_1,...,x_n]$, and if $K[\ul{x}]$ is an integral domain, $K(\ul{x}):=K(x_1,...,x_n)$ denotes the quotient field of $K[\ul{x}]$. For $\delta= (\delta_1,..., \delta_n) \in \mbb{Z}^n$, $\ul{x}^{\delta}$ denotes the element $x_1^{\delta_1}\cdot ... \cdot x_n^{\delta_n} \in K(\ul{x})$.

(ii) For any two subsets $A$, $B$ of a set $M$, we denote by $A \smallsetminus B$ the complement of $A \cap B$ in $A$. So $A \smallsetminus B = A \smallsetminus (A \cap B)$.  
\end{Nota}

\section{Kuroda's Example} \label{s2}

Let $K$ be a field of characteristic $0$, and let $K[x_1,...,x_4]=: K[\ul{x}]$ be the polynomial ring over $K$ in $4$ variables. 
Now for $i=1,2,3$, fix three elements $\delta_1, \delta_2, \delta_3 \in \mbb{Z}^4$ of the form  
\[\delta_1:=(-\delta_{1,1},\delta_{1,2},\delta_{1,3}, \delta_{1,4}), \quad \delta_2:=(\delta_{2,1},-\delta_{2,2},\delta_{2,3}, \delta_{2,4}), \quad \delta_3:=(\delta_{3,1},\delta_{3,2},-\delta_{3,3}, \delta_{3,4}), \]
where all $\delta_{i,4} \ge 0$, and all $\delta_{i,j} \ge 1$ for all $j \le 3$. Fix some $\gamma \in \mbb{Z}_{>0}$. Define for $i \in\{1,2,3\}$
\[ y_i := \underline{x}^{\delta_i}, \quad y_4 := x_4^{\gamma} \quad \text{and} \quad \Pi_i:= y_i - y_4,\] 
\[K[\underline{y}]:=K[y_1,...,y_4], \quad K[\underline{\Pi}]:=K[\Pi_1,\Pi_2,\Pi_3], \quad \text{and} \quad K(\underline{\Pi}):=K(\Pi_1,\Pi_2,\Pi_3),\]
\[R:= K[\ul{x}] \cap K(\ul{\Pi}).\] 
This notation differs somewhat from Kuroda's, in particular the variables named $y$ are not quite those that are named $y$ in \cite{MR2078390}. Using our notation we summarize some of the results from \cite{MR2078390}: 

\begin{Sum} \label{1} (Kuroda)
If the condition
\begin{equation} \label{q1}
\frac{\delta_{1,1}}{\delta_{1,1}+ min\{\delta_{2,1}, \delta_{3,1}\}}+ \frac{\delta_{2,2}}{\delta_{2,2}+ min\{\delta_{1,2}, \delta_{3,2}\}}
+ \frac{\delta_{3,3}}{\delta_{3,3}+ min\{\delta_{1,3}, \delta_{2,3}\}} < 1 
\end{equation}
holds, then the following obtains:

(i) The vectors $\delta_1, \delta_2, \delta_3, \gamma \mathbf{e}_4$ are linearly independent. Thus $K[\ul{y}]$ and $K[\ul{\Pi}]$ are polynomial rings, 
i.e. the generators are algebraically independent.

(ii) $R= K[\ul{x}] \cap K(\ul{\Pi}) = K[\ul{x}] \cap K[\ul{\Pi}]$

(iii) $R$ is not f.g. as a $K$-algebra, i.e. a counterexample to Hilbert's 14th problem. 
 
(iv) For all $i \neq j$ in $\{1,2,3\}$: $\delta_{i,i} \delta_{j,j} < \delta_{i,j} \delta_{j,i}$. $\quad$ \footnote{(iv) is actually not from
\cite{MR2078390}, but a (quite direct) consequence of equation (\ref{q1}).}    
\end{Sum}

\textbf{Concrete example:} One simple choice of the $\delta_i$ satisfying $(\ast)$ is: $\delta_1=(-1,3,3,0)$, $\delta_2=(3,-1,3,0)$, $\delta_3=(3,3,-1,0)$. 
When we talk about our \emph{concrete example} we will mean the one we get by choosing these $\delta_i$ and $\gamma =1$. 

\section{Semialgebraic realisation of $R$} \label{s3}

Let $R$ be any member of Kuroda's family of counterexamples to Hilbert's 14th problem over $K=\mbb{R}$, as described in Summary \ref{1}. We will construct an open semialgebraic $S \subset \mbb{R}^3$, such that $B_{\mbb{R}^3}(S)$ is isomorphic to $R$. 

\begin{Lemma} \label{2}
Let $V$, $V'$ be real affine $\mbb{R}$-varieties, $\Sigma' \subset V'$ semialgebraic. Let $f:V' \rightarrow V$ be a morphism induced by a homomorphism 
$\varphi: \mbb{R}[V] \rightarrow \mbb{R}[V']$ of the real coordinate rings, and set $\Sigma:=f(\Sigma')$. Then $\Sigma \subset V$ is semialgebraic, and:
\[B_V(\Sigma) = \varphi^{-1}\left(B_{V'}(\Sigma')\right)\]    
\end{Lemma}
   
\textbf{Proof:} That $\Sigma$ is semialgebraic follows from the Tarski-Seidenberg projection theorem (cf. \cite{MR1829790} Thm. 2.1.5.). 
$B_V(\Sigma) = \varphi^{-1}\left(B_{V'}(\Sigma')\right)$ is easy to check. \hfill $\square$

\textbf{Idea of construction of $S$:} We define the algebra $T:= \mbb{R}[\ul{y}] \cap \mbb{R}[\ul{x}]$. With this definition, obviously $R,T, \mbb{R}[\ul{\Pi}] \subset \mbb{R}[\ul{y}]$ and 
$R= T \cap \mbb{R}[\ul{\Pi}]$ (cf. Summary \ref{1} (ii)). One shows that $T$ is as $\mbb{R}$-algebra finitely generated by monomials in the variables $y_i$. This fact makes it easy to define an open semialgebraic $S' \subset \mbb{R}^4$ with $B_{\mbb{R}^4}(S')=T$, interpreting $\mbb{R}[\ul{y}]$ as the coordinate ring of $\mbb{R}^4$.
Now, by definition of $\mbb{R }[\ul{\Pi}]$, there is an inclusion $\varphi: \mbb{R}[\ul{\Pi}] \hookrightarrow \mbb{R}[\ul{y}], \; \Pi_i \mapsto y_i-y_4$, and if we interpret 
$\mbb{R}[\ul{\Pi}]$ as the coordinate ring of $\mbb{R}^3$, $\varphi$ induces a morphism $f: \mbb{R}^4 \rightarrow \mbb{R}^3$. We set $S:= f(S')$, 
and apply Lemma \ref{2} to get (inside $\mbb{R}[\ul{y}]$) :
\[B_{\mbb{R}^3}(S) = \varphi^{-1}(T) = T \cap \mbb{R}[\ul{\Pi}]= R.\]       
So $S$ is the example we wanted to construct. Because $f$ is just the projection parallel to the diagonal of $\mbb{R}^4$, we will be able tell how $S$ looks like.

The following Lemma provides the details missing in the idea of proof just given. 

\begin{Lemma} \label{3}
(i) $T$ is generated by the following set, which we denote by $M$
\[\left \{ \ul{y}^{\mathbf{n}} \; \middle | \; \mathbf{n}=(n_1,n_2,n_3,n_4) \in \mbb{Z}_{\ge 0}^4, \ \text{s.th.} \; \forall \{i,j,k\} = \{1,2,3\}: \; \delta_{i,i} n_i \le \delta_{j,i} n_j+ \delta_{k,i} n_k \right \} \]

(ii) We define 
\[ S':= \left \{(y_1,...,y_4) \in \mbb{R}^4 \; \middle | \;  \; \forall i,j \in \{1,2,3\} \; \text{with} \; i \neq j: \; |y_i^{\delta_{j,i}} y_j^{\delta_{i,i}}| < 1 , \; |y_4| < 1  \right\}.\]
With this definition:
\[\forall \lambda \in \mbb{R}_{>0} \quad B_{\mbb{R}^4}(\lambda S') = T, \quad \footnote{Where $\lambda S':= \{ x \in \mbb{R}^4 \mid \frac{1}{\lambda} x \in S'\}$} \]
as subalgebras of $\mbb{R}[\ul{y}]$.
\end{Lemma}

With these details filled in, we have proven that indeed $B_{\mbb{R}^3}(S)= R$, and even $B_{\mbb{R}^3}(\lambda S)= R$ for any $\lambda \in \mbb{R}_{>0}$. In particular $S$ is an open semialgebraic subset of $\mbb{R}^3$
such that $B_{\mbb{R}^3}(S)$ is not finitely generated as $\mbb{R}$-algebra.  
 
\textbf{Proof (Lemma \ref{3}):} (i): By definition of the $y_i$, a monomial in $\mbb{R}[\ul{y}]$ is contained in $\mbb{R}[\ul{x}]$ exactly if it is an element of $M$. A polynomial in $\mbb{R}[\ul{y}]$ 
lies in $\mbb{R}[\ul{x}]$ if and only if all monomials it is composed of lie in $\mbb{R}[\ul{x}]$.

(ii): First we prove $T \subseteq B_{\mbb{R}^4}(\lambda S')$ by showing that every $\ul{y}^{(n_1,...,n_4)} \in M$ is bounded on $\lambda S'$, for any $\lambda \in \mbb{R}_{>0}$. 
Take $(b_1,b_2,b_3,b_4) \in S'$ and WLOG assume $|b_1| > 1$ and thus $|b_2|, |b_3| < 1$. By definition of $S'$ we know 
\begin{equation} \label{e1}
 |b_1^{\frac{\delta_{2,1}}{\delta_{1,1}}} b_2| < 1, \quad |b_1^{\frac{\delta_{3,1}}{\delta_{1,1}}} b_3| < 1, \quad \text{and} \quad |b_4| < 1
\end{equation}
And by definition of $M$, $n_1 \le  \frac{\delta_{2,1}}{\delta_{1,1}} n_2 + \frac{\delta_{3,1}}{\delta_{1,1}} n_3$. Thus, with $n:=n_1+n_2+n_3+n_4$:
\[ |(\lambda b_1)^{n_1} (\lambda b_2)^{n_2} (\lambda b_3)^{n_3} (\lambda b_4)^{n_4}| \le \lambda^n |b_1^{n_1} b_2^{n_2} b_3^{n_3}| \le \lambda^n \left|b_1^{\left(\frac{\delta_{2,1}}{\delta_{1,1}} n_2 + \frac{\delta_{3,1}}{\delta_{1,1}} n_3 \right)} b_2^{n_2} b_3^{n_3}\right|\]   
\[= \lambda^n \left|\left(b_1^{\frac{\delta_{2,1}}{\delta_{1,1}}} b_2\right)^{n_2} \left(b_1^{\frac{\delta_{3,1}}{\delta_{1,1}}} b_3\right)^{n_3}\right| \le \lambda^n.\]
So $|\ul{y}^{(n_1,...n_4)}|$ is bounded by $\lambda^n$ on $\lambda S'$.

To prove that $B_{\mbb{R}^4}(\lambda S') \subseteq T$, consider any $g \in \mbb{R}[\ul{y}] \smallsetminus T$. Let $F$ be the set of all monomials $g$ consists of. WLOG we assume that for some of these monomials $\delta_{1,1} n_1  \le \delta_{2,1} n_2 + \delta_{3,1} n_3$ fails to hold.   
Among the monomials in $F$ select the monomial $h=\ul{y}^{(n_1,...,n_4)}$, for which the quadruple 
$\bigl( \delta_{1,1} n_1  - \delta_{2,1} n_2 - \delta_{3,1} n_3, \, -n_2, \, -n_3, \, -n_4 \bigr)$ is larger then for any other monomial in $F$, according to the lexicographical order.
(For this $h$, one has $\delta_{1,1} n_1  - \delta_{2,1} n_2 - \delta_{3,1} n_3 > 0$ by our WLOG assumption.) 

Using $\llog:= \log \circ \log$ and $\lllog:= \log \circ \log \circ \log$, consider the sequence 
\[a_k := (k^{\delta_{1,1}}, 1/(k^{\delta_{2,1}} \log k), 1/(k^{\delta_{3,1}} \llog k), 1/\lllog k).\] 
It lies inside $S'$, as we show below. Then,   
\[h(a_k)= k^{(\delta_{1,1} n_1  - \delta_{2,1} n_2 - \delta_{3,1} n_3)} (\log k)^{-n_2} (\llog k)^{-n_3} (\lllog k)^{-n_4} \xrightarrow{k \rightarrow \infty} \infty.\]
Furthermore, $h(a_k)$ goes to infinity at least by the order $\lllog k$ faster than $\widetilde{h}(a_k)$ for any $\widetilde{h} \in F \smallsetminus \{ h \}$. Hence $g(a_k) \rightarrow \infty$, and so $g \notin B_{\mbb{R}^4}( \lambda S')$. Thus $B_{\mbb{R}^4}(\lambda S') \subseteq T$.

It remains to show that $a_k \in \lambda S'$ (for $k \gg 0$). 
Let $a_{k,i}$ be the $i$-th component of the vector $a_k \in \mbb{Z}^4$. For $k \gg 0$, $|\frac{1}{\lambda} a_{k,4}| < 1$.
For $i\neq j \in \{1,2,3\}$ and $k \gg 0$, $|(\frac{1}{\lambda} a_{k,i})^{\delta_{j,i}} (\frac{1}{\lambda} a_{k,j})^{\delta_{i,i}}| < 1$ is obvious if none of $i,j$ is $1$. 
So it remains to check for all $j \in \{2,3\}$ that (for $k \gg 0$)
\[\left|\left(\frac{1}{\lambda} a_{k,1}\right)^{\delta_{j,1}} \left(\frac{1}{\lambda} a_{k,j}\right)^{\delta_{1,1}}\right| < 1 \quad \text{and} \quad \left|\left(\frac{1}{\lambda} a_{k,j}\right)^{\delta_{1,j}} \left(\frac{1}{\lambda} a_{k,1}\right)^{\delta_{j,j}}\right| 
< 1\] 
The first one is clear by definition of $a_{k,1}$ and $a_{k,2}$. For the second one, use 
$|(\frac{1}{\lambda} a_{k,j})^{\delta_{1,j}} (\frac{1}{\lambda} a_{k,1})^{\delta_{j,j}}| < (\frac{k}{\lambda})^{(\delta_{1,1} \delta_{j,j}- \delta_{1,j} \delta_{j,1})}$ and Summary \ref{1} (iv). \hfill $\square$

\textbf{How do $S$ and $S'$ look like?} In the $3$ dimensional subspace of $\mbb{R}^4$ defined by $y_4=0$ let $S''$ be the subset defined by   
\[S'':= \left \{(y_1,...,y_3) \in \mbb{R}^3 \; \middle | \; \forall i,j \in \{1,2,3\} \; \text{with} \; i \neq j: \; |y_i^{\delta_{j,i}} y_j^{\delta_{i,i}}| < 1   \right\}\]
$S''$ looks like a kind of star, whose center is at $(0,0,0)$ and whose $6$ infinitely long ``rays'' are lying on the coordinate axes.
Inside $\mbb{R}^4$, $S'$ is then the direct sum $S'' \oplus I$, where $I=\{(0,0,0,y_4) \, | \, |y_4| < 1 \}$. Since $f: \mbb{R}^4 \rightarrow \mbb{R}^3$ is the projection
along the diagonal $y_1=y_2=y_3=y_4$, $S= f (S')$ is the (not direct) sum $S = S'' + J$, where 
$J=\{(a,a,a) \in \mbb{R}^3 \, | \, |a| < 1 \}$. 

If we take $R$ to be the concrete example we fixed at the end of section \ref{s2}, then the star $S''= \left \{(y_1,...,y_3) \in \mbb{R}^3 \; \middle | \; \forall i,j \in \{1,2,3\} \; \text{with} \; i \neq j: \; |y_i^3 y_j| < 1 \right \}$ is quite symmetrical (cf. Figure 1).
\begin{figure}[H]
\centering
\includegraphics[width=4cm]{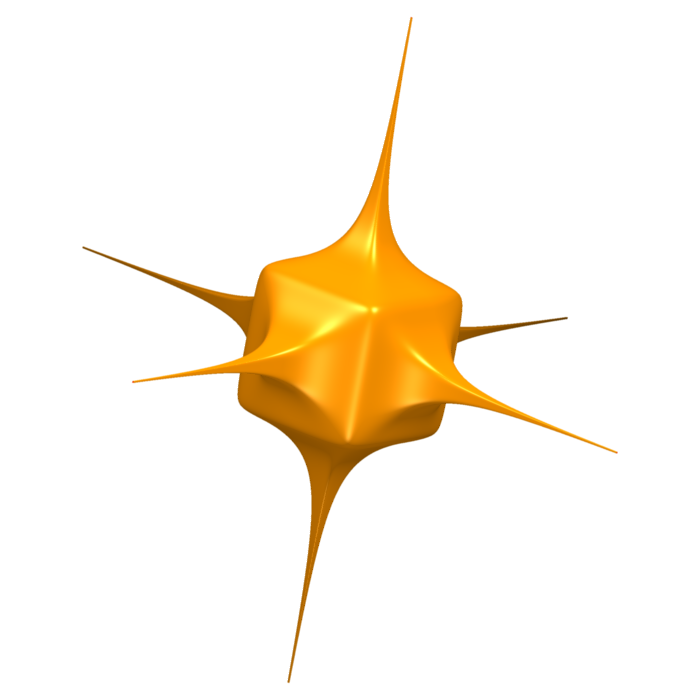}
$\qquad \qquad$ 
\includegraphics[width=4cm]{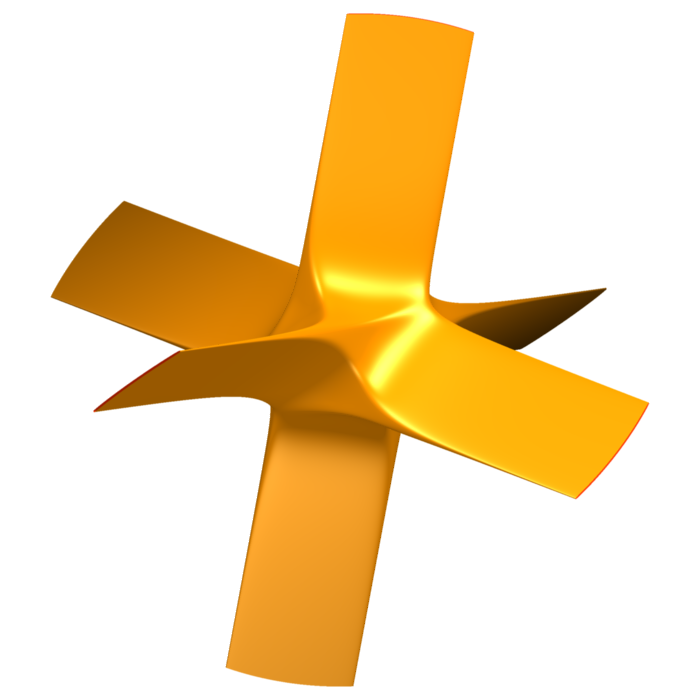}
\caption[]{Aproximate pictures of $S''$ resp. $\widetilde{S}$ for our concrete example \footnotemark }
\label{fig: a}
\end{figure}
\footnotetext{These pictures are actually images of algebraic sufaces (intersected with a ball around the origin) which approximate the boundaries of $S''$ resp. $\widetilde{S}$. The coordinate axes are aligned in the same way in both pictures, but the scale is not the same: $\widetilde{S}$ is viewed from a somewhat greater distance. The images where made using the free software SURFER (http://www.imaginary2008.de/surfer).}
We still can ``improve'' the semialgebraic set $S$ somewhat: It is not clear to me whether $S$ is basic, 
but we will now define an obviously basic open set $\widetilde S \subset \mbb{R}^3$, having the same ring of bounded polynomials. 

\begin{Def}\label{ds} (i) We set $d_i:= \min_{j \in \{1,2,3\} \smallsetminus \{i\}} \; \delta_{j,i}$. 

(ii) Let $\widetilde S \subset \mbb{R}^3$ be the basic open
set defined by the following inequalities.
\begin{gather}
(\Pi_1^{2 d_1} -1)(\Pi_2-\Pi_3)^{2 \delta_{1,1}} < 1 \tag{A1} \\
(\Pi_1^{2} -1)((\Pi_2+\Pi_3)^{2}-4) < 4 \tag{B1} \\
(\Pi_2^{2 d_2} -1)(\Pi_1-\Pi_3)^{2 \delta_{2,2}} < 1 \tag{A2} \\
(\Pi_2^{2} -1)((\Pi_1+\Pi_3)^{2}-4) < 4 \tag{B2} \\
(\Pi_3^{2 d_3} -1)(\Pi_1-\Pi_2)^{2 \delta_{3,3}} < 1 \tag{A3} \\
(\Pi_3^{2} -1)((\Pi_1+\Pi_2)^{2}-4) < 4 \tag{B3}
\end{gather}
\end{Def}

\begin{Prop} \label{pp}
Set $C:= \{(a_1,a_2,a_3) \in \mbb{R}^3 \; | \; |a_1| > 2, \; \text{or } \; |a_2| > 2, \; \text{or } \; |a_3| > 2 \}$,  then:

(i) $(\frac{1}{2} S) \cap C \subset \widetilde S \cap C \subset (2 S) \cap C$

(ii) Thus: $B_{\mbb{R}^3}(\widetilde S) = B_{\mbb{R}^3}(S)=R$
\end{Prop}

\textbf{Proof:} (i): Set 
\[S'':= \left \{(\Pi_1,...,\Pi_3) \in \mbb{R}^3 \; \middle | \; \forall i,j \in \{1,2,3\} \; \text{with} \; i \neq j: \; |\Pi_i^{\delta_{j,i}} \Pi_j^{\delta_{i,i}}| < 1 \right \}.\]
Then $S=S''+J$ where $J:= \{(a,a,a) \in \mbb{R}^3 \; | \; |a| < 1 \}$, as explained at the end of Section \ref{s3}. 

Take any $(a_1,a_2,a_3) \in \frac{1}{2} S \cap C$, WLOG $|a_1| > 2$. We can write $(a_1,a_2,a_3)=(b_1+a,b_2+a,b_3+a)$ 
with $|a| < \frac{1}{2}$, $(b_1,b_2,b_3) \in \frac{1}{2} S''$. Since $|b_1| > 1$ this implies $|b_2|, |b_3| < \frac{1}{2}$, and thus 
$|a_2+a_3|=|b_2+b_3+2a| < 2$. Hence $(a_1,a_2,a_3)$ satisfies inequality (B1). We get (A1) by:
\[ |(2b_1)^{\delta_{2,1}} (2 b_2)^{\delta_{1,1}}| < 1 \; \text{and } \;  |(2b_1)^{\delta_{3,1}} (2 b_3)^{\delta_{1,1}}| < 1\]
\[\Rightarrow |(2 b_1)^{d_1} (2 \max \{b_2,b_3 \})^{\delta_{1,1}}| < 1 \; \Rightarrow |(a_1)^{d_1}(a_2-a_3)^{\delta_{1,1}}| < 1\]
\[\Rightarrow (a_1)^{2d_1}(a_2-a_3)^{2\delta_{1,1}} < 1 \; \Rightarrow  ((a_1)^{2d_1}-1) (a_2-a_3)^{2\delta_{1,1}} < 1 \]
That the remaining (A2)-(B3) are fulfilled is easy to check. We have shown the first inclusion of (i). 

Now take any $(a_1,a_2,a_3) \in \widetilde{S} \cap C$, WLOG $|a_1| > 2$ and WLOG $\delta_{2,1}=\min\{\delta_{2,1}, \delta_{3,1}\} = d_1$. Then (A1) resp. (B1) imply $|a_2+a_3| < \frac{5}{2}$, $|a_2-a_3| < 1$ hence $|a_3| < 2$, so $(a_3,a_3,a_3) \in 2J$. Thus it suffices to show $(b_1,b_2,0):=(a_1-a_3,a_2-a_3,a_3-a_3) \in 2 S''$.   
We only have to check the two defining inequalities of $2 S''$ not containing $\Pi_3$. 
Under our two WLOG assumptions:                           
\[ (A1) \ \Rightarrow \ |(a_1^{2d_1}-1)b_2^{2\delta_{1,1}}| < 1 \ \Rightarrow \ |a_1^{2 \delta_{2,1}}b_2^{2\delta_{1,1}}| < 2\]
\[ \Rightarrow \  |a_1^{\delta_{2,1}}b_2^{\delta_{1,1}}| < 2 \ \Rightarrow \ |(\frac{1}{2} b_1)^{\delta_{2,1}}b_2^{\delta_{1,1}}| < 2 \ \Rightarrow \ |(\frac{1}{2} b_1)^{\delta_{2,1}}(\frac{1}{2} b_2)^{\delta_{1,1}}| < 1 \]
This is one of the defining inequalities not containing $\Pi_3$, and together with Summary \ref{1} (iv) and $|b_2|=|a_2-a_3|<1$ it also implies the second one by:
\[1 >  |(\frac{1}{2}b_1)^{\delta_{2,1}} (\frac{1}{2}b_2)^{\delta_{1,1}}|^{\frac{\delta_{1,2}}{\delta_{1,1}}} = |(\frac{1}{2}b_1)^{\delta_{2,1}\frac{\delta_{1,2}}{\delta_{1,1}}} (\frac{1}{2}b_2)^{\delta_{1,2}}| \ \Rightarrow \ |(\frac{1}{2}b_1)^{\delta_{2,2}} (\frac{1}{2}b_2)^{\delta_{1,2}}|<1 \]
 
(ii): Part (i) implies that $B_{\mbb{R}^3}(2 S) \subseteq B_{\mbb{R}^3}(\widetilde S) \subseteq B_{\mbb{R}^3}(\frac{1}{2} S)$. But we have seen above that $B_{\mbb{R}^3}(\lambda S)= R$ for any $\lambda \in \mathbb{R}_{>0}$. \hfill $\square$

\begin{Rem}
This easily generalises to $B_{\mbb{R}^3}(\lambda \widetilde S) = B_{\mbb{R}^3}(S) =R$ for all $\lambda \in \mbb{R}_{>0}$. Using this one can check (with some effort) that the right hand sides of the inequalities defining $\widetilde{S}$, may be replaced by any real numbers for (B1), (B2), (B3), and by any positive real numbers for (A1), (A2), (A3), without changing $B_{\mbb{R}^3}(\widetilde{S})=R$. If one of the latter however is replaced by a non-positive number, $B_{\mbb{R}^3}(\widetilde{S})$ will become strictly larger than $R$.
\end{Rem}

\section{Quasiprojective realisation of $R$}

Over any field $K$ of characteristic $0$, we construct, by blowing up a $\mbb{P}^3_K$, a projective variety $X$ containing $\mbb{A}^3_K$, such that
there is an open subvariety $\mbb{A}^3_K \subset U \subset X$ 
with $\GO(U) =R$. In the special case $K = \mbb{R}$ this construction produces a completion $\mbb{A}^3_{\mbb{R}} \hookrightarrow X$,
which is compatible with $\widetilde S$ as well as $S$. (The idea for the construction of $X$ stems from this special case.)

\begin{Lemma} \label{4} Any   $f(\Pi_1, \Pi_2, \Pi_3) \in K[\ul{\Pi}]$ can of course also be expressed in the form 
$f=f^{[1]}(\Pi_1, \Pi_2, \Pi_2-\Pi_3)$, $f=f^{[2]}(\Pi_2, \Pi_1, \Pi_1-\Pi_3)$
and $f=f^{[3]}(\Pi_3, \Pi_1, \Pi_1 - \Pi_2)$. We can write 
\[f= f^{[1]}(\Pi_1, \Pi_2, \Pi_2-\Pi_3)= \sum_{(r_1,r_2,r_3) \in \mbb{Z}_{\ge 0}^3} \alpha_{(r_1,r_2,r_3)} \Pi_1^{r_1} \Pi_2^{r_2} (\Pi_2-\Pi_3)^{r_3}\]
and we set $m_1(f):=\{(r_1,r_2,r_3) \in \mbb{Z}_{\ge 0}^3 \; | \;  \alpha_{(r_1,r_2,r_3)} \neq 0 \}$. 
We define $m_2(f)$ resp. $m_3(f)$ analogously, using the decomposition of $f^{[2]}$ resp. $f^{[3]}$ into monomials. 
Then we can describe $R \subset K[\ul{\Pi}]$ as the set of all $f \in K[\ul{\Pi}]$ fulfilling the following condition:
\[ \forall i \in \{1,2,3\} \quad \bigl( (r_1,r_2,r_3) \in m_i(f) \; \Rightarrow \; \delta_{i,i} r_1 \le d_i r_3 \bigr) \tag{$\ast$} \]
\end{Lemma} 

\textbf{Proof:}  First consider any $f \in K[\ul{\Pi}] \smallsetminus R$. Then $f$ WLOG contains, when expressed in the variables $y_i$, 
a monomial $\ul{y}^{(n_1,...,n_4)}$ such
that $\delta_{i,i} n_1 > \delta_{2,1} n_2 + \delta_{3,1} n_3$ (cf. Lemma \ref{3} (i)). Now there must be a $(r_1,r_2,r_3) \in m_1(f)$ such that the monomial $\ul{y}^{(n_1,...,n_4)}$ appears in 
$\Pi_1^{r_1} \Pi_2^{r_2} (\Pi_2-\Pi_3)^{r_3}$, i.e. in $(y_1-y_4)^{r_1} (y_2-y_4)^{r_2} (y_2-y_3)^{r_3}$.
This is only possible if $r_1 \ge n_1$ and $r_3 \le n_2+n_3$. Hence:
\[ \delta_{1,1} r_1 \ge \delta_{1,1} n_1 > \delta_{2,1} n_2 + \delta_{3,1} n_3 \ge d_1 (n_2 +n_3) \ge d_1 r_3 \]  
Thus condition ($\ast$) is violated, as it should be. 
Now consider any $f$ violating ($\ast$). WLOG there is a $(r_1,r_2, r_3) \in m_1(f)$, satisfying $\delta_{1,1} r_1 > d_1 r_3$.  
Among these, choose the $(r_1,r_2, r_3)$ which is maximal according to the lexicographical order. Then, in    
$(y_1-y_4)^{r_1} (y_2-y_4)^{r_2} (y_2-y_3)^{r_3}$,
the terms $y_1^{r_1} y_4^{r_2} y_2^{r_3}$ and $y_1^{r_1} y_4^{r_2} y_3^{r_3}$ appear. It is easy to check, considering how $(r_1,r_2, r_3)$ was chosen, that these
two monomials are not canceled out in $f$, and that at least one of them violates the inequality $\delta_{1,1} n_1 \le \delta_{2,1} n_2 + \delta_{3,1} n_3$. Thus $f \notin R$. \hfill $\square$ 

Lemma \ref{4} provides conditions on an $f \in K[\ul{\Pi}]$ for being an element of $R$. 
We will construct the completion $\mbb{A}_{K}^3 \hookrightarrow X$ in such a way that those 
conditions are equivalent to $f$ not having poles along certain components of the boundary $Y:=X \smallsetminus \mbb{A}_{K}^3$.  
Rewrite the inequalities in the conditions as $\frac{r_1}{r_3} \le \frac{d_i}{\delta_{i,i}}$. If $\frac{d_i}{\delta_{i,i}} \in \mbb{N}$ for all $i \in \{1,2,3\}$, it is more easy to see how to construct such an $X$ by blowups of $\mbb{P}_K^3$. In general $\frac{d_i}{\delta_{i,i}}$ is not an integer. In that case the idea how to still construct a suitable completion $X$ is basically to use Euclid's algorithm.

\textbf{Notation:} For $Q$ a rational number, we denote by $\langle Q \rangle$ the fractional part of $Q$, i.e. $\langle Q \rangle := Q - \lfloor Q \rfloor$.    

For $i \in \{1,2,3\}$, set $Q_i:= \frac{d_i}{\delta_{i,i}}$. For each $i$ we define numbers $q_{i,m}$ (compare to Euclid's algorithm) by  
\[q_{i,1}:= \lfloor Q_i \rfloor, \; q_{i,2}:=  \lfloor \langle Q_i \rangle^{-1} \rfloor, \; q_{i,3}:=  \lfloor \langle \langle Q_i \rangle^{-1} \rangle^{-1} 
\rfloor, \; ..., \; q_{i,M_i}:=  \lfloor \langle \langle  ... \langle Q_i \underbrace{\rangle^{-1} ...\rangle^{-1} \rangle^{-1}}_{(M_i-1)-times} \rfloor \]
where $M_i$ is the (first) number such that 
\[q_{i,M_i}:=  \lfloor \langle \langle  ... \langle Q_i \rangle^{-1} ...\rangle^{-1} \rangle^{-1} \rfloor = 
\langle \langle  ... \langle Q_i \rangle^{-1} ...\rangle^{-1} \rangle^{-1}.\] 

For $N_i:= \sum_{m=1}^{M_i} q_{i,m}$, divide the set $\{0,1,...,N_i\}$ into the sets
\[I_{i,1}:=\{0,...,q_{i,1}\}, \; I_{i,2}:=\{q_{i,1} +1,...,q_{i,1}+q_{i,2} \}, \; ..., \; I_{i,M_i}:= \{(\sum_{m=1}^{M_i-1} q_{i,m}) +1,..., N_i\}.\]

\textbf{Construction of the completion $X$:} Consider $K[\Pi_1,\Pi_2, \Pi_3]$ as the coordinate ring of $\mbb{A}^3_K$, embed $\mbb{A}^3_K$ in the projective space
$\mbb{P}^3_K$ with homogeneous coordinates $(z_1:z_2:z_3:z_4)$ such that $z_i/z_4= \Pi_i$. Set 
\[B:= \mbb{P}^3_K \smallsetminus \mbb{A}^3_K, \quad p_1:= (1:0:0:0), \ p_2:= (0:1:0:0), \ p_3:= (0:0:1:0) \in \mbb{P}^{3}_K.\] 
We construct $X$ from $\mbb{P}_K^3$, by preforming at each of the $3$ points $p_i$ a series of $N_i+1$ blowups, which will be defined below. We enumerate these blowups starting with $0$, and we call the  exceptional divisor introduced by the $n-th$ blowup $E_{i,n}$. In the $0$-th step, we blow up the point $p_i$. This introduces the exceptional divisor $E_{i,0}$. After this and also after every following step of our series of blowups, we will denote the strict transforms of $B$ and of the exceptional divisors $E_{i,k}$ again by the same symbols. 
To avoid having to distinguish different cases when describing the blowups, we introduce $I_{i,0}:=\{-1\}$ and $E_{i,-1}:= \mc{V}_{\mbb{P}^3_K}(z_k-z_l)$, where $\{k,l\} = \{1,2,3\} \smallsetminus \{i\}$. For any $n \in \{-1,0,...,N_i\}$ let $k(n)$ be the number such that $n \in I_{i,k(n)}$.
Now we define recursively, for $n \in \{0,1,..., N_i-1\}$, the $n+1$-th blowup of the series of blowups around $p_i$: 

Set $l_{n}:= \max I_{i,k(n)-1}$, then the center of the $n+1$-th blowup is the (reduced) intersection $C_{n+1}:= E_{i,n} \cap E_{i,l_{n}}$.

We denote the boundary $X \smallsetminus \mbb{A}_K^3$ by $Y$. Set 
\[J_{1,i}:= \{0,1,..., N_i-1\}, \qquad J_{2,i}:= \bigcup_{m \in \{1,...,M_i\}, \text{$m$ odd}} I_{i,m} \smallsetminus \{N_i\},\]
\[Z_1:= B \cup \bigcup_{i=1}^{3} \bigcup_{n \in J_{1,i}} E_{i,n}, \qquad Z_2:= B \cup \bigcup_{i=1}^{3} \bigcup_{n \in J_{2,i}} E_{i,n}.\]  

\begin{Prop} \label{5}

(i) For $K$ any field of characteristic $0$, as subsets of $K[\ul{\Pi}]$,
\[\GO(X \smallsetminus Z_1) = \GO(X \smallsetminus Z_2) = R. \]
In particular the rings $\GO(X \smallsetminus Z_1) = \GO(X \smallsetminus Z_2)$ are not f.g. over $K$.

(ii) For $K = \mbb{R}$, the inclusion of $\mbb{A}^3_K$ in $X$ is a completion compatible with $S$ as well as $\widetilde S$ (in the sense of Definition \ref{d1} (ix)).
More precisely the intersection of the closure of $S$ resp. $\widetilde S$ in $X$ with $Z_2$ is empty, and the intersection with the components of $Y$ not contained in $Z_2$ is Zariski-dense in those components. 
 
\end{Prop}

\textbf{Proof:} (i): For any $f \in K(X)$ let $\Pole(f) \subset X$ be the set of poles of $f$. It is of pure codimension $1$. Using the description of $R$ in Lemma \ref{4}, it thus suffices to prove that for all $i \in \{1,2,3\}$ the following conditions 
on an $f \in K[\ul{\Pi}] \subset K(X)$ are equivalent.

($\Cond_{i,1}$) $(r_1,r_2,r_3) \in m_i(f) \; \Rightarrow \; \delta_{i,i} r_1 \le d_i r_3$

($\Cond_{i,2}$) For all the exceptional divisors $E_{i,n}$: $E_{i,n} \subseteq \Pole(f) \; \Rightarrow \; E_{i,n} \subseteq Z_1$.

($\Cond_{i,3}$) For all the exceptional divisors $E_{i,n}$: $E_{i,n} \subseteq \Pole(f) \; \Rightarrow \; E_{i,n} \subseteq Z_2$.

We only show this in the case $i=1$, the others being analogous. 
At first we show the equivalence only for $f$ of the form $f = \Pi_1^{r_1} \Pi_2^{r_2} (\Pi_2-\Pi_3)^{r_3}$. Let 
\[X_{N_1} \overset{\varphi_{N_1}}{\longrightarrow} X_{N_1-1} \overset{\varphi_{N_1-1}}{\longrightarrow} ....\overset{\varphi_1}{\longrightarrow}  X_0 \overset{\varphi_0}{\longrightarrow} \mbb{P}^3_K\]
denote the series of $N_1+1$ blowups preformed at the point $p_1$, as introduced above in defining $X$. \footnote{We ignore the blowups at the points $p_2$ and $p_3$, since they are irrelevant for evaluating the three condictions for $i=1$.}         
We will iteratively define charts $A_n \subset X_n$ with $A_n \cong \mbb{A}_K^n$ together with coordinates on each $A_n$. These charts and coordinates will enable us to check for each $E_{1,n}$ whether  $E_{1,n} \subseteq \Pole(f)$.   
First take the affine chart $A \cong \mbb{A}_K^3$ of $\mbb{P}^3_K$ defined by $z_1 \neq 0$. We choose coordinates 
$\alpha:= z_2/z_1$, $\beta:= (z_2-z_3)/z_1$, $\gamma:= z_4/z_1$ on $A$, and with this choice $p_1= (0,0,0)$. 
Now $f$ can be considered as the rational function $z_1^{r_1} z_2^{r_2} (z_2-z_3)^{r_3} z_4^{-(r_1+r_2+r_3)}$  on $\mbb{P}^3_K$. On $A$ we can express this as 
$f(\alpha,\beta, \gamma)=\alpha^{r_2} \beta^{r_3} \gamma^{-(r_1+r_2+r_3)}$. 

Equip $A \times \mbb{P}^2_K$ with coordinates $(\alpha, \beta, \gamma, a:b:c)$. The transform $\widetilde A:= \varphi_0^{-1}(A)$ can be seen as the subvariety of $A \times \mbb{P}^2_K$ which is defined by the equations $\alpha b = \beta a$, $\alpha c = \gamma a$, $\beta c= \gamma b$. 
Let $A_{0}$ be the subvariety of $\widetilde A$ defined by $c \neq 0$. Let $\alpha_0,\beta_0, \gamma_0$ be coordinates of $\mbb{A}_K^3$ 
and define an open embedding $\mbb{A}_K^3 \hookrightarrow \widetilde A$ by $(\alpha_0, \beta_0, \gamma_0) \mapsto (\gamma_0 \alpha_0, \gamma_0 \beta_0, \gamma_0, \alpha_0: \beta_0: 1)$. 
This is an isomorphism with $A_0$, and using it, we consider $(\alpha_0, \beta_0, \gamma_0)$ as coordinates of $A_0$. When we pull back $f(\alpha, \beta, \gamma)$ to $A_0$ it becomes
$f_0(\alpha_1,\beta_0,\gamma_0):=f(\gamma_0 \alpha_0, \gamma_0 \beta_0, \gamma_0)= \alpha_0^{r_2} \beta_0^{r_3} \gamma_0^{-r_1}$.     


For the exceptional divisor $E_{1,0}= \varphi_0^{-1}(p_1)$, the subvariety $E_{1,0} \cap A_0$ of $A_0$ (which we call $E_{1,0}$ again) is defined by $\gamma_0 =0$, i.e. $E_{1,0}= \mc{V}_{A_0}(\gamma_0)$. Also $E_{1,-1}= \mc{V}_{A_0}(\beta_0)$.\footnote{Where $E_{1,-1}$ is the strict transform of the original $E_{1,-1}$, restricted to $A_0$.}

Now we will define recursively the further affine charts $A_n$. The ``$n \rightarrow n+1$'' step of the recursion (for $n < N_1$), we start with a given $\mbb{A}_K^3 \cong A_n \subset X_n$ with coordinates $(\alpha_n, \beta_n, \gamma_n)$, satisfying the following ``\emph{induction hypothesis}'': Depending on whether $k(n)$ is odd (case $(1)$), or even (case $(2)$):
\[(1): \; E_{1,n}=\mc{V}_{A_n}(\gamma_n), \; E_{1,l_{n}} = \mc{V}_{A_n}(\beta_n), \quad \text{or},\]
\[ \quad (2): \; E_{1,n}=\mc{V}_{A_n}(\beta_n), \; E_{1,l_{n}} = \mc{V}_{A_n}(\gamma_n).\]
For $n=0$ we have defined $A_0$ and coordinates $(\alpha_0, \beta_0, \gamma_0)$, satisfying this hypothesis. 

Hence, if we denote the restriction to $A_n$ of the center $C_{n+1}$ of the next blowup again by $C_{n+1}$, then $C_{n+1}= \mc{V}_{A_n}(\beta_n,\gamma_n)$, and  $\widetilde A_n:= \varphi_{n+1}^{-1}(A_n)$ is isomorphic to the blowup of $A_n$ in $C_{n+1}$. Choosing coordinates $(\alpha_n, \beta_n, \gamma_n, b:c)$ on $A_n \times \mbb{P}^1_K$, we can consider $\widetilde{A}_n$ as the subvariety described by $\beta_n c = \gamma_n b$. The new exceptional divisor is $E_{1,n+1}= \mc{V}_{\widetilde A_n}(\beta_n,\gamma_n)$. And, depending on the case :
\[(1): \; E_{1,n}=\mc{V}_{\widetilde A_n}(c), \; E_{1,l_{n}}=\mc{V}_{\widetilde A_n}(b), \quad \text{or},\]
\[ \quad (2): \; E_{1,n}=\mc{V}_{\widetilde A_n}(b), 
\; E_{1,l_{n}}=\mc{V}_{\widetilde A_n}(c).\] 

We distinguish two further cases now: Either $k(n+1) = k(n)$ (case $(A)$), or $k(n+1)= k(n)+1$ (case $(B)$). Together with the distinction between (1) and (2), these combine to $4$ possible cases $(1A)$, $(1B)$, $(2A)$, $(2B)$. We call $(O)$ the union of the cases $(1A)$, $(2B)$, and $(E)$ the union of $(1B)$, $(2A)$, since 
$k(n+1)$ is odd iff $(O)$, even iff $(E)$.

We want our next affine chart $A_{n+1} \subset \widetilde{A}_{n} \subset X_{n+1}$ to contain open dense subsets of $E_{n+1}$ and $E_{l_{n+1}}$ 
(and of $C_{n+2} =E_{n+1} \cap E_{l_{n+1}}$ if $n+2 \le N_1$). 
In case $(A)$, $l_{n+1} = l_{n}$, In case $(B)$, $l_{n+1} = n$. 
We choose $A_{n+1}:= \widetilde{A}_n \smallsetminus \mc{V}_{\widetilde{A}_n}(c)$ in case $(O)$, and $A_{n+1}:= \widetilde{A}_n \smallsetminus \mc{V}_{\widetilde{A}_n}(b)$ in case $(E)$.

Let $(\alpha_{n+1}, \beta_{n+1},\gamma_{n+1})$ be coordinates of $\mbb{A}_K^3$. In case $(O)$ the embedding $\mbb{A}_K^3 \hookrightarrow \widetilde A_n$ defined by
$(\alpha_{n+1}, \beta_{n+1}, \gamma_{n+1}) \mapsto (\alpha_{n+1}, \beta_{n+1} \gamma_{n+1}, \gamma_{n+1}, \beta_{n+1}:1)$, is an isomorphism with $A_{n+1}$, and we use it to 
give $A_{n+1}$ the coordinates $(\alpha_{n+1}, \beta_{n+1}, \gamma_{n+1})$. 
In case $(E)$ the same holds for the embedding defined by 
$(\alpha_{n+1}, \beta_{n+1}, \gamma_{n+1}) \mapsto ( \alpha_{n+1}, \beta_{n+1}, \beta_{n+1} \gamma_{n+1}, 1: \gamma_{n+1})$.

One checks that, in these new coordinates,
depending on case $(O)$ or $(E)$:  
\[(O): \; E_{1,n+1}=\mc{V}_{A_{n+1}}(\gamma_{n+1}), \; E_{1,l_{n+1}} = \mc{V}_{A_{n+1}}(\beta_{n+1}), \quad \text{or},\]
\[(E): \; E_{1,n+1}=\mc{V}_{A_{n+1}}(\beta_{n+1}), \; E_{1, l_{n+1}} = \mc{V}_{A_{n+1}}(\gamma_{n+1})\]
I.e. the induction hypothesis is fulfilled for $A_{n+1}$ with coordinates $(\alpha_{n+1}, \beta_{n+1}, \gamma_{n+1})$. 
 
Now let $f_n$ be the pullback to $A_n$ of the rational function $f$ on $\mbb{P}^3$. 
We already pulled $f$ back to $A_0$ and obtained $f_0= \alpha_0^{r_2} \beta_0^{r_3} \gamma_0^{-r_1}$. Using our inductive description of the blowups just given, 
for $n \in \{0,...,N_1-1\}$ we obtain the following pullback formulas:
\[ f_{n+1}(\alpha_{n+1}, \, \beta_{n+1}, \, \gamma_{n+1}) =  f_n(\alpha_{n+1}, \, \beta_{n+1} \gamma_{n+1}, \, \gamma_{n+1}), \quad \text{if $k(n+1)$ odd}, \tag{$\dagger$}\]
\[ f_{n+1}(\alpha_{n+1}, \, \beta_{n+1}, \, \gamma_{n+1}) =  f_n(\alpha_{n+1}, \, \beta_{n+1}, \, \beta_{n+1} \gamma_{n+1}), \quad \text{if $k(n+1)$ even.} \tag{$\dagger$}\] 

Let $r_1^{(n)}$, $r_2^{(n)}$, $r_3^{(n)}$ be the integers such that $f_{n}= \alpha_{n}^{r^{(n)}_2} \beta_{n}^{r^{(n)}_3} \gamma_{n}^{-r^{(n)}_1}$, in particular $r_i^{(0)} = r_i$. Since on the chart $A_n$, an open part of $E_{1,n}$ coincides with an open part of $\mc{V}_{A_n}(\gamma_n)$ or $\mc{V}_{A_n}(\beta_n)$ (depending on the parity of $k(n)$), we have: If $k(n)$ is odd, $E_{i,n} \subseteq \Pole(f)$ iff $r_1^{(n)} > 0$, while if $k(n)$ is even,
$E_{i,n} \subseteq \Pole(f)$ iff $r_3^{(n)} < 0$. 

Set $\nu(m):= \max I_{1,m}$. The pullback formulas ($\dagger$) yield,
\[\text{if $m$ odd}: \; \left(r_1^{(\nu(m))}, \, r_2^{(\nu(m))}, \, r_3^{(\nu(m))} \right)= \left(r_1^{(\nu(m-1))}- q_{1,m} r_3^{(\nu(m-1))}, \, r_2, \, r_3^{(\nu(m-1))} \right)\]
\[\text{if $m$ even}: \; \left(r_1^{(\nu(m))}, \, r_2^{(\nu(m))}, \, r_3^{(\nu(m))} \right)= \left(r_1^{(\nu(m-1))}, \, r_2, \, r_3^{(\nu(m-1))}-q_{1,m} r_1^{(\nu(m-1))}\right)\]
   
Now we get, with $Q_1$ still denoting $\frac{d_1}{\delta_{1,1}}$: 
\[\Cond_{1,1} \; \Leftrightarrow \; r_1 \le Q_1 r_3 \; \Leftrightarrow \; r_1-q_{1,1}r_3 \le Q_1 r_3-q_{1,1} r_3= \langle Q_1 \rangle r_3 \; 
\Leftrightarrow \; r_1^{(\nu(1))} \le \langle Q_1 \rangle r_3^{(\nu(1))}\]
\[ \Leftrightarrow \; \langle Q_1 \rangle^{-1} r_1^{(\nu(1))}  \le r_3^{(\nu(1))} \; \Leftrightarrow \; 
\langle \langle Q_1 \rangle^{-1} \rangle r_1^{(\nu(2))}  \le r_3^{(\nu(2))}\]
\[\Leftrightarrow \; ... \; \Leftrightarrow \; r_1^{(\nu(3))} \le \langle \langle \langle Q_1 \rangle^{-1} \rangle^{-1} \rangle r_3^{(\nu(3))} \]  
In case $M_1$ is odd the chain of equivalences continues until
\[... \; \Leftrightarrow \; r_1^{(\nu(M_1-1))} - q_{1,M_1} r_3^{(\nu(M_1-1))} \le \langle \langle  ... \langle Q_1 \underbrace{\rangle^{-1} ...\rangle^{-1} \rangle^{-1}}_{(M_1-1)-times} r_3^{(\nu(M_1-1))}-q_{1,M_1} r_3^{(\nu(M_1-1))} \] 
\[ \Leftrightarrow \; r_1^{(N_1)} \le 0 \; \Leftrightarrow \Cond_{1,2}\]
In case $M_1$ is even, the chain ends instead with: $... \ 0 \le r_3^{(N_1)} \; \Leftrightarrow \; \Cond_{1,2} $. 

Using parts of this chain of eqivalences, it is also easy to show recursively that $r_3^{(n)} \ge 0$ for all $n \le N_1$, given $\Cond_{1,1}$.
Hence $\Cond_{1,1} \Rightarrow \Cond_{1,3}$. The proof of the equivalence for $f = \Pi_1^{r_1} \Pi_2^{r_2} (\Pi_2-\Pi_3)^{r_3}$ is finished, as  $\Cond_{1,3} \Rightarrow \Cond_{1,2}$ is trivial.  

Now for a general $g \in K[\ul{\Pi}]$,
$g= \sum_{\nu=1}^{r} s_{\nu} f_{\nu}$ for some $s_{\nu} \in K$ and suitable $f_{\nu} = \Pi_1^{r_{1,\nu}} \Pi_2^{r_{2,\nu}} (\Pi_2-\Pi_3)^{r_{3,\nu}}$, and 
we claim that each $\Cond_{1,i}$ holds for $g$ iff it holds for all the $f_{\nu}$. For $\Cond_{1,1}$ this is obvious. For $\Cond_{1,2}$, $\Cond_{1,3}$, it follows from the fact that $g_n= \sum_{\nu=1}^{r} s_{\nu} \alpha_{n}^{r^{(n)}_{2, \nu}} \beta_{n}^{r^{(n)}_{3,\nu}} \gamma_{n}^{-r^{(n)}_{1,\nu}}$ and that 
$(r_{1,\nu}^{(n)},r_{2,\nu}^{(n)}, r_{3,\nu}^{(n)}) = (r_{1,\nu'}^{(n)},r_{2,\nu'}^{(n)}, r_{3,\nu'}^{(n)})$ implies $(r_{1,\nu},r_{2,\nu}, r_{3,\nu}) = (r_{1,\nu'},r_{2,\nu'}, r_{3,\nu'})$.    

(ii): We show that for any $\lambda \in \mbb{R}_{>0}$, the completion $\mbb{A}_{\mbb{R}}^3 \hookrightarrow X$ is compatible with $\lambda \widetilde{S}$, 
and that the closure $\lambda \widetilde{S}^{cl}$ of $\lambda \widetilde{S}$ meets exactly the boundary divisors not contained in $Z_2$. This implies that the same is true for $S$, by Proposition \ref{pp} (i). Obviously $\lambda \widetilde{S}^{cl}$ does not meet $B$. We now treat the divisors $E_{1,n}$, the other $E_{i,n}$ being analogous.

Now  $\lambda \widetilde{S}$ is defined by inequalities $( \lambda A1)$, $(\lambda B1)$,..., analogous to those in Definition \ref{ds} (ii). Expressed in coordinates $(\alpha, \beta,\gamma)$, $( \lambda A1)$ reads:
\[\lambda^{-2(d_1+\delta_{1,1})} \beta^{2 \delta_{1,1}} \gamma^{-2 (d_1+ \delta_{1,1})}- \lambda^{-2\delta_{1,1}} \beta^{2 \delta_{1,1}} \gamma^{-2 \delta_{1,1}} < 1.\]
\[\text{Pulled back to $A_0$:} \quad  \lambda^{-2(d_1+\delta_{1,1})} \beta_0^{2 \delta_{1,1}}\gamma_0^{-2 d_1} - \lambda^{-2\delta_{1,1}} \beta_0^{2 \delta_{1,1}} < 1. \tag{$\ddagger$}\]
Using ($\dagger$), we obtain that the pullback of the left had side of ($\ddagger$) has a pole along any $E_{1,n} \subset Z_2$ in $X$. Thus $\lambda \widetilde{S}^{cl}$ does not meet any $E_{1,n} \subset Z_2$. It remains to show that it meets Zariski-dense any $E_{1,n} \nsubseteq Z_2$. We instead show this for the smaller set $\widehat S \subset \lambda \widetilde S$ we define next: 

Note that the conjunction of $(\Pi_2/\lambda)^2 < 1/4$ and 
$((\Pi_2-\Pi_3)/\lambda)^2 < 1/4$, implies also $(\Pi_3/\lambda)^2 < 1$ and so implies all defining inequalities of $\lambda \widetilde S$, except of (A1). So let $\widehat S$ be the subset of $\mbb{R}^3$ defined by $(\lambda A1)$, $(\Pi_2/\lambda)^2 < 1/4$ and $((\Pi_2-\Pi_3)/\lambda)^2 < 1/4$. 

Extended to $\mbb{P}^3$ and then pulled back to $A_0$, the latter two conditions become $\alpha_0^2 < \frac{\lambda^2}{4}$ and $\beta_0^2 < \frac{\lambda^2}{4}$. Using ($\dagger$), we see that on any divisor $E_{1,n}$ for $n \ge 1$ the pullback of the function $\beta_0^2$ is zero. Furthermore $\alpha_0^2 < \lambda^2$ pulls back to $\alpha_n^2 < \lambda^2$. 
Also we obtain that on $(Y \smallsetminus Z_2) \smallsetminus E_{1,N_1}$, the pullback of the left hand side of ($\ddagger$) is zero. 
Thus, the closure of $\widehat S$ intersects each component in the Zariski-dense subset defined by $\alpha_n^2 < \lambda^2$. Finally on $E_{1,{N_1}}$,  
depending on whether $k(N_i)$ is odd or even, one gets from ($\ddagger$) instead a condition of the form $\lambda^{-2(d_1+\delta_{1,1})} \beta_{N_1}^s < 1$ resp. 
$\lambda^{-2(d_1+\delta_{1,1})} \gamma_{N_1}^{-s} < 1$ for some $s >0$. This together with $\alpha_{N_1} < \lambda^2$ still defines a Zariski-dense subset. \hfill $\square$ 

\textbf{Remark:} For the concrete example specified at the end of section \ref{s2}, we have $M_1=M_2=M_3=1$ and $N_1=N_2=N_3=3$, and thus $Z_1=Z_2$.  

\section{Appendix}
            
\subsection*{Characterisations of rings of bounded polynomials on regular semialgebraic subsets of normal varieties} 

The following Proposition generalizes a weaker version of Theorem 4.5. in \cite{PS1} (cf. Summary \ref{su} (iii))
from dimension $\le 2$ to arbitrary dimension. The proof is very similar to the proof in \cite{PS1}, we just do allow the embedded resolution of singularities
to also blow up parts of the original variety $V$.       

\begin{Prop} \label{ab}
Let $V$ be a real normal irreducible quasiprojective $\mbb{R}$-variety, let $S \subset V$ be semialgebraic and regular. 
Then there is an open subvariety $V' \subseteq V$, such that, $B_V(S) = B_{V'}(S')$ for $S' := S \cap V'$, and such that $V'$ has a
completion compatible with $S'$. In particular there is a nonsingular quasiprojective $\mbb{R}$-variety $U$, birational to $V$, such that $\GO(U) \cong B_V(S)$, as $\mbb{R}$-algebras. \footnote{The last sentence of the proposition also applies if $S$ is only regular \emph{at infinity}, since then $B_V(S)= B_V(\inter(S))$.}  
\end{Prop}

\textbf{Proof:} If we speak of desingularisation of a variety in the following, we will always mean resolution of singularities by a series of 
blowups in regular closed centers, with properties as for example described in \cite{MR1978567} (Chapter -1). 
Set $S_{reg}:=S \cap V_{reg}$. Since $V$ is normal $\GO (V_{reg})=\GO (V)$ and since in addition $S$ is regular, 
$B_{V_{reg}}(S_{reg})=B_V(S)$. Let $V_{reg} \hookrightarrow X$ be a completion of $V_{reg}$, with $X$ nonsingular, such that $B:= X \smallsetminus V_{reg}$ 
is a divisor. Let $S^{cl}$ be the closure of $S_{reg}$ in $X(\mbb{R})$, let $\partial S^{cl}$ be its boundary and let $C$ be the Zariski-closure of $\partial S^{cl}$in $X$.
Since $S$ is regular, $C$ is a divisor in $X$. Let $Z$ be the set of singularities of $C$. 
Set $V':= V_{reg} \smallsetminus (V_{reg} \cap Z)$, $S':= S_{reg} \cap V'$. Again $B_{V'}(S') =B_V(S)$, since $V$ is normal and 
only a codimension $\ge 2$ subset of $V$ is missing in $V'$. 

Now let $\varphi: \widetilde X \rightarrow X$ be the embedded desingularisation of $B \cup C$ inside $X$. On $V'$, $\varphi$ is an isomorphism, and we denote $\varphi^{-1}(V')$ by $V'$ again. Set $\widetilde B:= \widetilde X \smallsetminus V'$, $D:= \varphi^{-1} (B \cup C)$, then $\widetilde B \subseteq D$. The closure 
of $\widetilde S:= \varphi^{-1} (S')$ in $\widetilde X (\mbb{R})$ is $\widetilde S^{cl}=\varphi^{-1}(S^{cl})$.  
 
We claim that the completion $V' \hookrightarrow \widetilde X$ is compatible with $S'$, which also implies the rest of our Proposition via part (ii) of Summary 
\ref{su}. We have $\widetilde B =\varphi^{-1} (B \cup Z)= \widehat{B} \cup W$, where $\widehat{B}$ is the strict transform of $B$, $W$ the exceptional divisor of $\varphi$.  So $\widetilde{B}$ is a normal crossings divisor, by the properties of desingularisation and the fact that $B$ is a divisor.
Hence to show compatibility of the completion, it suffices to check that for any irreducible component $E$ of $\widetilde B$, 
$\widetilde S^{cl} \cap E$ is either Zariski-dense in $E$ or empty. 
Choose some $p \in \widetilde S^{cl} \cap E$, if the set is nonempty. Also $D= \varphi^{-1}(B \cup C)$ is a normal crossings divisor. 
So a neighborhood of $p$ in $\widetilde X(\mbb{R})$ looks like some ball in $\mbb{R}^n$, and the components of $\varphi^{-1}(B \cup C)$ which contain $p$ (among them $E$) meet like coordinate hyperplanes in $p$. They cut the ball into sectors adjacent to $p$. 
If one restricts to a small enough neighborhood of $p$, the regular $\widetilde S^{cl}$ has to contain one of these sectors, since $\varphi^{-1}(C)$ contains the boundary of $\widetilde S^{cl}$. But the boundary of each sector contains a Zariski-dense part of $E$, so $\widetilde S^{cl} \cap E$ is Zariski-dense in $E$. \footnote{If you want a stricter formulation of the argument of this last paragraph cf. the proof of Thm. 4.5. in \cite{PS1}.} \hfill $\square$
   
It is known that for a $K$-algebra $R$, the following conditions are equivalent. \footnote{In the proof of these equivalences one can use a more technical equivalent condition (VI), whose role it is, to imply condition (I) and to be implied by all other conditions. It is: (VI) $R$ is a Krullring, isomorphic to the subalgebra of a f.g. integral $K$-algebra, and for every prime ideal $\mathfrak{p}$ of $R$ of hight $1$, the quotient field of $R/\mathfrak{p}$ has transcendence degree over $K$ exactly one less then the quotient field of $R$. }  
(where $\cong$ is isomorphism of $K$-algebras) 

(I) There is a normal irreducible quasiaffine $K$-Variety $U$ such that $\GO(U) \cong R$.

(II) There is a normal irreducible $K$-Variety $U$ such that $\GO(U) \cong R$.

(III) For some field extension $K|L|E$, and some normal f.g. $K$-algebra $A$ contained in $E$, $R \cong L \cap A$.

(IV) $R \cong  A_1 \cap ... \cap A_n$ where the $A_i$ are normal f.g. $K$-algebras contained in some extension field of $K$. ($n \in \mbb{N}$) 

(V) $R$ is normal and integral, and there are elements $f_1,...,f_n \in R$, such that the localisations $R_{f_i}$ all are f.g. $K$-algebras, and such that
$R = R_{f_1} \cap ... \cap R_{f_n}$.

The equivalence of $(I)$ up to $(IV)$ can be gathered together from work of M. Nagata in the context of Hilbert's 14th problem (cf. \cite{MR0088034} and \cite{MR0096645}). 
J. Winkelmann reproved most of these equivalences in \cite{MR1981881} (Thm. 1 and Thm. 2), and this may be a more compact source for the results.  
$(V)$ is equivalent to the others, since obviously $(V) \Rightarrow (IV)$, and since $(I) \Rightarrow (V)$.
   
Proposition \ref{ab} implies, that for $V$ an irreducible normal $\mbb{R}$-variety and for $S \subset V$ semialgebraic and regular, the $\mbb{R}$-algebra $B_V(S)$, always has the $5$ equivalent properties just listed (with $K= \mbb{R}$).  In fact the following stronger characterisation holds.

\begin{Char} For an $\mbb{R}$-algebra $R$ the following conditions are equivalent:

(i) $R$ is isomorphic to the algebra $B_V(S)$ for some irreducible normal affine real $\mbb{R}$-variety $V$ and some regular semialgebraic $S \subset V$.

(ii) The equivalent properties (I)-(V) listed above hold, and in addition the quotient field of $R$ is (formally) real.

(The same holds with ``regular'' replaced by ``open'' or ``regular at infinity'' in (i).)
\end{Char}
      
\textbf{Proof:}  The ``only if'' direction is just what was explained above, together with the fact that $\mbb{R}(V)$ of a real variety $V$ is (formally) real. 
For the ``if'' direction, first choose, by property (I), a normal irreducible quasiaffine $K$-Variety $U$ such that $\GO(U) \cong R$,  
which has to be real since $\mbb{R}(U)$ is the quotient field of $\GO(U)=R$ and thus (formally) real. Then apply Thm 4.11. of \cite{PS1}. \hfill $\square$   

\begin{small}
\bibliographystyle{alpha}
\addcontentsline{toc}{chapter}{References}
\bibliography{BibSeb}
\end{small}

\end{document}